\documentclass{article}

\usepackage{endnotes}
\usepackage{amssymb}
\usepackage{german}
\usepackage{natbib} 

\usepackage{times, mathptm}

\makeatletter
\def\enoteheading{\section*{\notesname
  \@mkboth{\uppercase{\notesname}}{\uppercase{\notesname}}}%
     \leavevmode\par
}
\def\enoteformat{\rightskip\z@ \leftskip\z@ \parindent=\z@
     \leavevmode{\hbox{{\@theenmark}.\quad}}}
\makeatother

\title{Hilbert's `\emph{Verungl\"uckter Beweis}', the first epsilon
theorem, and consistency proofs} \author{Richard Zach\\ Department of
Philosophy\\ University of Calgary\\ 2500 University Drive N.W.\\
Calgary, Alberta T2N 1N4, Canada\\ rzach@ucalgary.ca}

\def\Bcite#1#2{\citet[#1]{#2}}
\def\pp#1{#1}

\let\citeN\citet
\let\citeA\citealt
\let\cite\citep
\let\frak\mathfrak
\let\note\endnote
\let\ol\overline

\def\ft{\frak{t}}
\def\fA{\frak{A}}\def\fB{\frak{B}}

\def\sg{\selectlanguage{german}}
\def\eg{\selectlanguage{USenglish}}
\eg
\let\epsilon\varepsilon

\begin{document}
\bibliographystyle{hpl}

\maketitle

\begin{abstract}
In the 1920s, Ackermann and von Neumann, in pursuit of Hilbert's
Programme, were working on consistency proofs for arithmetical
systems.  One proposed method of giving such proofs is Hilbert's
epsilon-substitution method. There was, however, a second approach
which was not reflected in the publications of the Hilbert school in
the 1920s, and which is a direct precursor of Hilbert's first epsilon
theorem and a certain `general consistency result' due to Bernays.  An
analysis of the form of this so-called `failed proof' sheds further
light on an interpretation of Hilbert's Programme as an instrumentalist
enterprise with the aim of showing that whenever a `real' proposition
can be proved by `ideal' means, it can also be proved by `real',
finitary means.
\end{abstract}

\section{Introduction}

The aim of Hilbert's Programme for consistency proofs in the 1920s is
well known: to formalize mathematics, and to give finitary
consistency proofs of these systems and thus to put mathematics on a
`secure foundation'. What is perhaps less well known is exactly how
Hilbert thought this should be carried out.  Over ten years before
Gentzen developed sequent calculus formalizations of arithmetic and
used an elaboration of his cut-elimination procedure to give a
consistency proof of Peano Arithmetic, Hilbert proposed a different
approach. He believed that the principles criticized by intuitionists,
the principle of the excluded middle in its application to infinite
totalities and the use of unbounded existential quantifiers are, at
root, the same. This root is the axiom of choice.  In a course on the
foundations of mathematics, he remarked that whereas the use of
unbounded quantification results in significant problems for giving a
consistency proof,
\begin{quotation}
the core of the difficulty lies at a different point, to which one
usually only pays attention later: it lies with Zermelo's \emph{axiom
of choice}\dots We want to extend the axiom of choice.  To each
proposition with a variable $A(a)$ we assign an object for which the
proposition holds only if is holds in general.  So, a counterexample,
if one exists.\note{\citeN[30--31]{Hilbert:22b}. Following
\citeN{HilbertBernays:34}, we will use the following notation: $a$,
$b$, \dots{} stand for free variables, whereas $x$, $y$, \dots are
bound variables.  $A$, $B$, \dots{} are formula variables. $\frak{A}$,
$\frak{B}$, \dots are metavariables for formulas,
and $\frak{n}$, $\frak{z}$ stand for numerical terms.  For
uniformity, the notation in some quotations has been adjusted.}
\end{quotation}
This counterexample is given by the $\tau$-operator: $\tau_x A(x)$ is
an object for which $A(x)$ is false, if there is one. The dual
operator $\epsilon_x A(x)$, is a witness, i.e., an object for which
$A(x)$ is true, if $A(x)$ is true for anything.\note{The
$\tau$-operator was mentioned in \citeN{Hilbert:22a} and formally
introduced, together with the transfinite axiom, in
\citeN{Hilbert:23}.  In a course given in Winter 1922/23
\cite{Hilbert:22b,Hilbert:22c}, Hilbert and Bernays introduce the
$\epsilon$-operator for the first time, although initially with the
same interpretation (as a counterexample) as the $\tau$-operator.  In
notes appended to the typescript \citeyearpar{Hilbert:22b}, the final
form of the $\epsilon$-operator and $\epsilon$-axioms appears for the
first time.  Kneser's notes to the course \cite{Hilbert:22c} do not
contain this final version, suggesting that this change was made after
the conclusion of the 1922/23 course.}  The $\epsilon$-operator is
governed by the \emph{transfinite axiom},
\[
A(a) \to A(\epsilon_x A(x)).
\]
A finitary consistency proof of mathematical theorems which allows
the elimination of applications of the choice principle (in the form
given to it by the transfinite axiom) would then show that such
application is justified after all.  It would also show that unbounded
quantification is admissible in mathematics, since with the help of
the transfinite axioms one can define quantifiers by
\[
(\exists x) A(x) \equiv A(\epsilon_x A(x)) \quad \textrm{and}\quad 
(\forall x) A(x) \equiv A(\epsilon_x \neg A(x)).
\]
$\epsilon$-terms may be seen as ideal elements whose addition to the
theory of finite propositions reintroduces the powerful methods of
infinite mathematics, and `round out the theory'.  To show that
their addition is permissible requires a proof that $\epsilon$-terms
can be eliminated from proofs of `real', finitary, propositions.
This elimination of $\epsilon$-terms from formal proofs in
arithmetical theories was to proceed according to the
epsilon-substitution method.  Hilbert's approach here was to define a
finitary procedure which would produce, given a proof involving
$\epsilon$-terms, a substitution of these terms by actual
numbers.\note{The basic idea was presented in \citeN{Hilbert:23} and
in the course mentioned \cite{Hilbert:22b,Hilbert:22c}, for
discussion, see \citeN{Zach:03}.  Roughly, the idea is this: first
replace every $\epsilon$-term by~0. The instances of the transfinite
axiom for an $\epsilon$-term $\epsilon_x \frak{A}(x)$ in the proof
then become formulas of the form $\frak{A}(\frak{n}) \to \frak{A}(0)$.
If this formula is false, $\frak{A}(\frak{n})$ is true.  In the next
iteration of the procedure, replace $\epsilon_x \frak{A}(x)$ by
$\frak{n}$.  The difficulty is to extend this idea to the case where
more than one $\epsilon$-term, and in particular, when nested
$\epsilon$-terms occur in the proof.}  Applying this substitution to
the proof would then result in a purely elementary proof about numbers
which would contain no trace of the transfinite elements of the
original proof.  In addition, it is seen finitarily that all
initial formulas, and hence also the end formula, of the resulting
proof are true. Since such a proof cannot possibly have a
contradiction as its last line, the consistency of arithmetic would be
established.  Hilbert presented his `Ansatz' for finding such
substitutions in \citeN{Hilbert:22a}; it was extended by
\citeN{Ackermann:24} and \citeN{Neumann:27}.

The epsilon-substitution method and its role in Hilbert's Programme are
now relatively well understood.  There was, however, a \emph{second}
proposal for proving consistency, also based on the epsilon calculus,
which has escaped historical attention, and which was never presented
in the publications of the Hilbert school before 1939.  In the second
volume of \emph{Grundlagen der Mathematik} \cite{HilbertBernays:39},
Bernays first developed in print a well worked-out theory of the
epsilon calculus as an alternative formulation and extension of
predicate logic, and proved the so-called first and second epsilon
theorems.  In Section 1.4, Bernays presented a `general consistency
theorem' based on the first epsilon theorem, which applies, e.g., to
elementary geometry and to arithmetic with an open induction rule.
This second approach to consistency proofs via the first epsilon
theorem, however, dates back to the beginning of Hilbert's
Programme.  In a letter from Bernays to Ackermann of October 1929,
Bernays refers to this second approach as Hilbert's
`\emph{verungl\"uckter Beweis}' (the `failed proof').  This failed
proof is not mentioned in Hilbert's publications of the early
1920s nor into his lectures on the subject of 1922 and 1923.  A record
of the basic idea, a second `\emph{Ansatz}', is, however, available
in the form of a six-page note in Bernays's hand.

The aim of this paper is to present and analyze this second approach
to proving consistency, and to show how Hilbert's `verungl\"uckter
Beweis' precipitated the later proof of the first epsilon theorem by
Bernays and Ackermann.  Given the role envisaged by Hilbert for the
$\epsilon$-calculus and the consistency proofs based on it, such an
analysis will help illuminate not just the genesis of an important
proof-theoretic result (the epsilon theorem), but also Hilbert's aim
and strategy for providing consistency proofs.  In the following
section, we will revisit the first epsilon theorem, and show how it
can be used to establish consistency results. Following this
discussion, I present the suggestion contained in Hilbert's second
\emph{Ansatz}, and outline why this approach was not pursued by
Hilbert and his students in the 1920s.  A concluding section discusses
the relevance of the result in the wider context of Hilbert's
consistency project.

\section{The first epsilon theorem and the general consistency result}

The epsilon calculus consists in the elementary calculus of free
variables plus the `transfinite axiom', $A(a) \to A(\epsilon_x A(x))$.
The elementary calculus of free variables is the quantifier-free
fragment of the predicate calculus, i.e., axioms for propositional
logic and identity, with substitution rules for free individual ($a$,
$b$, \dots) and formula ($A$, $B$, \dots) variables and modus ponens.

One of the most basic and fruitful results concerning Hilbert's
$\epsilon$-calculus is the so-called epsilon theorem.  It states that
if a formula $\frak{E}$ containing no $\epsilon$-terms is derivable in
the $\epsilon$-calculus from a set of axioms which also do not contain
$\epsilon$-terms, then $\frak{E}$ is already derivable from these
axioms in the elementary calculus of free variables (i.e., essentially
using propositional logic alone).  A relatively easy consequence of
this theorem (or rather, of its proof) is Herbrand's theorem, and, in
fact, one of the first published correct proofs of Herbrand's
theorem is that given by Bernays in \emph{Grundlagen der Mathematik
II} \cite{HilbertBernays:39} based on the first
$\epsilon$-theorem.\note{\citeN{Gentzen:34} mentions that a version of
Herbrand's Theorem is a consequence of his `Versch\"arfter
Hauptsatz'. He does not, however, spell out the details.}
\citeN{Leisenring:69} even formulates the $\epsilon$-theorem in such a
way that the connection to Herbrand's theorem is obvious:
\begin{quotation}
If $E$ is a prenex formula derivable from a set of prenex formulas
$\Gamma$ in the predicate calculus, then a disjunction $B_1 \lor
\ldots \lor B_n$ of substitution instances of the matrix of $E$ is
derivable in the elementary calculus of free variables from a set
$\Gamma'$ of substitution instances of the matrices of the formulas
in~$\Gamma$.
\end{quotation}
Even without this important consequence, which was of course not
discovered until after Herbrand's \citeyearpar{Herbrand:30} thesis, the first
$\epsilon$-theorem constitutes an important contribution to
mathematical logic.  Without the semantical methods provided by the
completeness theorem for predicate logic, it is not at all clear that
the addition of quantifiers in the guise of $\epsilon$-terms and the
axioms governing them is a conservative extension of the elementary
calculus.  Keeping in mind the role of epsilon-terms as `ideal
elements' in a proof, the eliminability of which is the main aim of a
consistency proof of any mathematical system formulated with the aid
of the epsilon calculus, the first epsilon theorem is also the main
prerequisite for such a consistency proof.

Bernays stated the first and second epsilon theorem as follows:
\begin{quotation}
These theorems both concern a formalism $F$, which results from the
predicate calculus by adding to its symbols the $\epsilon$-symbol and
also certain individual [constant], predicate, and function symbols,
and to its axioms the $\epsilon$-formula [the transfinite axiom] and
furthermore certain \emph{proper axioms $\frak{P}_1, \ldots,
\frak{P}_\frak{k}$ which do not contain the $\epsilon$-symbol.}  For
such a formalism $F$, the two theorems state the following:
\begin{enumerate}
\item If $\frak{E}$ is a formula derivable in~$F$ which does not
contain any bound variables, and the axioms $\frak{P}_1, \ldots,
\frak{P}_\frak{k}$ also contain no bound variables, then the formula
$\frak{E}$ can be derived from the axioms $\frak{P}_1, \ldots,
\frak{P}_\frak{k}$ without the use of bound variables at all, i.e.,
with the elementary calculus of free variables alone (`first epsilon
theorem').
\item If $\frak{E}$ is a formula derivable in $F$ which does not
contain the $\epsilon$-symbol, then it can be derived from the axioms
$\frak{P}_1, \ldots, \frak{P}_\frak{k}$ without the use of the
$\epsilon$-symbol, i.e., with the predicate calculus alone (`second
epsilon theorem'). \cite[18]{HilbertBernays:39}
\end{enumerate}
\end{quotation}
The predicate calculus is formulated with a substitution rule for free
individual and formula variables; the elementary calculus of free
variables is the quantifier-free fragment of the predicate calculus
(i.e., without quantifier axioms or rules) or equivalently, the
epsilon calculus without transfinite axioms and without defining
axioms for the quantifiers.

A proof that the $\epsilon$-calculus is conservative over the
elementary calculus of free variables in the way specified by the
first epsilon theorem constitutes a proof of consistency of the
$\epsilon$-calculus and of mathematical theories which can be
formulated in such a way that the first $\epsilon$-theorem applies
(i.e., the axioms are quantifier- and $\epsilon$-free).  Indeed
Bernays used the first $\epsilon$-theorem to prove a `general
consistency result' for axiomatic theories to which the first
$\epsilon$-theorem applies.  Let us first outline the consistency
proof for a very basic arithmetical theory.  This theory results from
the elementary calculus of free variables by adding the constant $0$
and successor ($+1$) and predecessor ($\delta$) functions.  The
additional axioms are:
\begin{eqnarray*}
0 & \neq & x +1\\
x & = & \delta(x +1)
\end{eqnarray*}
To prove that the resulting axiom system is consistent, assume there
were a proof of $0 \neq 0$.  First, by copying parts of the
derivation as necessary, we can assume that every formula in the proof
is used as the premise of an inference at most once.  Hilbert and
Bernays call this `resolution into proof threads'.  The resulting
proof is in tree form; a branch of this tree (beginning with an axiom
and ending in the end-formula) is a \emph{proof thread}.  Next, we can
substitute numbers for the free variables in the proof (`elimination
of free variables').  Bernays describes this as follows:
\begin{quotation}
We follow each proof thread, starting at the end formula, until we
reach two successive formulas $\frak{A}$, $\frak{B}$ where the first
results from the second by substitution.  We record the substitution
also in the formula $\frak{B}$, so that we get instead of $\frak{B}$
a repetition of the formula~$\frak{A}$.

If $\frak{B}$ is an initial formula [axiom], then the substitution has been
transferred to the initial formula.  Otherwise, $\frak{B}$ was obtained by
substitution into a formula $\frak{C}$ or by repetition, or as conclusion of
an inference
\[
\begin{array}{cc}
\frak{C} & \quad \frak{C} \to \frak{B} \\
\multicolumn{2}{c}{\diagdown \quad \diagup}\\
\multicolumn{2}{c}{\frak{B}.}
\end{array}
\]
In the first case, we in turn replace $\frak{C}$ by $\frak{A}$, so
that the substitutions leading from $\frak{C}$ to $\frak{B}$ and from
$\frak{B}$ to $\frak{A}$ are recorded simultaneously. (In the case of
repetition, only one substitution is recorded.)

In the case of the inference schema [modus ponens], we record the
substitution leading from $\frak{B}$ to $\frak{A}$ in the formulas
$\frak{C}$ and $\frak{C} \to \frak{B}$; this changes the formula
$\frak{C}$ if and only if it contains the variables being substituted
for in the transition from $\frak{B}$ to $\frak{A}$.  In any case, the
original inference schema with conclusion $\frak{B}$ is replaced by an
inference schema
\[
\begin{array}{cc}
\frak{C}^* & \quad \frak{C}^* \to \frak{A} \\
\multicolumn{2}{c}{\diagdown \quad \diagup}\\
\multicolumn{2}{c}{\frak{A}.}
\end{array}
\]
We can proceed in this way until we reach an initial formula in each
thread.  When the procedure comes to its end, each substitution has
been replaced by a repetition, each inference schema by another
inference schema, and certain substitutions have been applied to the
initial formulas. \cite[225]{HilbertBernays:34}
\end{quotation}
Remaining free variables can now be replaced by $0$ (for individual
variables) and $0 = 0$ (for formula variables).  We would thus obtain
a proof of $0 \neq 0$ without free variables.

If we now reduce the variable-free terms in the resulting proofs to
standard numerals by successively replacing $\delta(0)$ by $0$ and
$\delta(\frak{t} + 1)$ by $\frak{t}$, we get a proof where each
initial formula is either an instance of a tautology, of an identity
axiom, or, if the original formula was one of the axioms for $+1$ and
$\delta$, a formula of the form of either
\begin{eqnarray*}
0 & \neq & \frak{n} + 1\\
\frak{n} & = & \frak{n}
\end{eqnarray*}
(where $\frak{n}$ is either $0$ or of the form $0+\cdots+1$).

Call an equation of the form $\frak{n} = \frak{n}$ `true' and one of
the form $\frak{n} = \frak{m}$, where $\frak{n}$ and $\frak{m}$ are
not identical, `false'. This can be extended to propositional
combinations of equations in the obvious way.  We observe that the
resulting proof has all true initial formulas, and since modus ponens
obviously preserves truth, all other formulas are also true. Since $0
\neq 0$ is false, there can be no proof of $0 \neq 0$.

Hilbert presented this proof in his course on the foundations of
mathematics in 1921/22 \citeyearpar{Hilbert:21,Hilbert:21a} and outlined
the basic approach in his 1922 talk at the meeting of the Deutsche
Naturforscher-Gesellschaft in Leipzig \citeyearpar{Hilbert:22c}. In a
course given the following year \citeyearpar{Hilbert:22b,Hilbert:22c},
Bernays and he extended it to axioms for primitive recursive
functions; \citeN{Ackermann:24} further elaborated it to include
second-order primitive recursive functions (see \citeA{Zach:03}).  The
challenge was to extend it to the case where $\epsilon$-terms and the
transfinite axiom are also present, leading to Hilbert's
$\epsilon$-substitution method. There, the aim was to find
substitutions not just for the free variables, but also for the
$\epsilon$-terms, ultimately also resulting in a proof without free or
bound variables and with true initial formulas.  An alternative method
is this: Instead of treating $\epsilon$-terms together with other
terms of the system, eliminate $\epsilon$-terms \emph{first}. We
introduce a step at the beginning of the proof which reduces a proof
in the $\epsilon$-calculus to one in the elementary calculus of free
variables as in the first $\epsilon$-theorem.  Thus, with the first
$\epsilon$-theorem in hand, Bernays could later formulate the following
`general consistency theorem':
\begin{quotation}
Let $F$ be a formalism which results from the predicate calculus by
adding certain individual, predicate, and function symbols.  Suppose
there is a method for determining the truth value of variable-free
atomic formulas uniquely. Suppose furthermore that the axioms do not
contain bound variables [i.e., no quantifiers and no $\epsilon$-terms]
and are verifiable [i.e., every substitution instance is true].  Then
the formalism is consistent in the strong sense that every derivable
variable-free formula is a true
formula. \cite[36]{HilbertBernays:39}
\end{quotation}
Suppose $\frak{A}$ is variable-free and derivable in~$F$.  If a
formalism $F$ satisfies the conditions, then the first
$\epsilon$-theorem yields a proof of $\frak{A}$ already in the
elementary calculus of free variables.  The procedures above
(resolution into proof threads, elimination of free variables) yields a
proof of $\frak{A}$ from substitution instances of the axioms of~$F$.
Since the axioms of $F$ are verifiable, these substitution instances
are true, and again, modus ponens preserves truth. So $\frak{A}$ is
true.  The requirement that the truth-value of variable-free atomic
formulas is decidable ensures that this is a finitary proof: It can
be finitarily verified that any \emph{given} proof in fact has
true initial formulas (and hence, a true end formula).

\section{Hilbert's \emph{Verungl\"uckter Beweis}}

The first $\epsilon$-theorem is first formulated in print in
\citeN{HilbertBernays:39}, but Hilbert had something like it in mind
already in the early/mid 1920s.  When working on \emph{Grundlagen der
Mathematik} in the late 1920s, Bernays revisited the idea, which had
been abandoned in favour of the $\epsilon$-substitution method.  In
correspondence with Ackermann in 1929 (discussed below), Bernays
refers to `Hilbert's second consistency proof for the
$\epsilon$-axiom, the so-called `failed proof', and suggests ways in
which the difficulties originally encountered could be fixed.
Surprisingly, this `failed proof', a precursor of the first
$\epsilon$-theorem, is not to be found in the otherwise highly
interesting elaborations of lecture courses on logic and proof theory
given by Hilbert (and Bernays) between 1917 and 1923.  The only
evidence that the $\epsilon$-theorem predates Bernays's proof of it in
\citeN{HilbertBernays:39} are the letter from Bernays to Ackermann from
1929, and a sketch of the simplest case of the theorem.

The sketch in question is a six-page manuscript in Bernays's hand
which can be found bound with the lecture notes to Hilbert's course
\emph{Elemente und Prinzipienfragen der Mathematik}, taught in the
Summer Semester 1910 in G\"ottingen \cite{Hilbert:10}.  Although it
bears a note by Hilbert `Insertion in WS [Winter Semester] 1920 [sic]', the
notation used in the sketch suggests that it was written after
sometime after 1922, when the epsilon notation was first
introduced. Unfortunately, the only substantial discussion of the
proof is found in the letter from Bernays to Ackermann from 1929
quoted below.  The fact that it uses the $\epsilon$-axioms in their
final form, suggests that it was written after Hilbert and Bernays's
1922/23 course, in which $\epsilon$s were still used in their dual
forms. One would also expect Hilbert to have presented the proof in
the course, if it had had been available then.  The proof is briefly
alluded to in a letter from Ackermann to Bernays of June 1925 in a way
that suggests that it was not a recent proof.\note{Ackermann to
Bernays, June 25, 1925. 9pp.  ETH Z\"urich Library, Hs. 975.96.
Ackermann writes: `[Critical formulas] are eliminated in the way in
which Hilbert wanted to earlier [\emph{fr\"uher}]' (p.~3). Ackermann
here tries to use the idea of the `failed proof '
to fix his own faulty $\epsilon$-substitution proof.}  Thus, the proof
likely dates from 1923 or 1924.

The sketch bears the title `Consistency proof for the logical axiom
of choice $Ab \to A\epsilon_a Aa$, in the simplest case'.  In it, we
find a proof of the first $\epsilon$-theorem for the case where the
substitution instances of the transfinite axiom used in the proof,
i.e., the so-called \emph{critical formulas}
\[
\frak{A}(\frak{t}) \to \frak{A}(\epsilon_x \frak{A}(x))
\]
are such that $\frak{A}(x)$ contains no $\epsilon$'s, and furthermore
the identity axioms are not used at all.  The proof goes as follows.
Suppose
\begin{eqnarray*}
\frak{A}(\frak{t_1}) & \to & \frak{A}(\epsilon_x \frak{A}(x))\\
 & \vdots \\
\frak{A}(\frak{t_n}) & \to & \frak{A}(\epsilon_x \frak{A}(x))\\
\end{eqnarray*}
are all the critical formulas involving~$\frak{A}$ in a proof
of~$0 \neq 0$. First, replace every formula $\frak{F}$ occurring in
the proof by the conditional $\ol{\fA}(\ft_1) \to \frak{F}$, and every
application of modus ponens by the (derivable) inference
\[
{\ol\fA(\ft_1) \to \frak{S} \qquad \ol\fA(\ft_1) \to (\frak{S} \to
\frak{T})} \over \ol\fA(\ft_1) \to \frak{T}\] Every formula resulting
thus from a substitution instance~$\frak{F}$ of an axiom (other than
the critical formula for $\frak{t}_1$) is then derivable by
\[
{\frak{F} \qquad \frak{F} \to (\ol\fA(\ft_1) \to \frak{F})}
\over 
{\ol\fA(\ft_1) \to \frak{F}}
\]
The formula corresponding to the $\epsilon$-axiom involving $\ft_1$ is
derived using
\[
\begin{array}{c}
\fA(\ft_1) \to (\ol\fA(\ft_1) \to \fA(\epsilon_x \fA(x)) \\
(\fA(\ft_1) \to (\ol\fA(\ft_1) \to \fA(\epsilon_x \fA(x))) \to 
(\ol\fA(\ft_1) \to (\fA(\ft_1) \to \fA(\epsilon_x \fA(x))) \\
\hline
\ol\fA(\ft_1) \to (\fA(\ft_1) \to \fA(\epsilon_x \fA(x))
\end{array}
\]
The premises of this inference are propositional axioms.  Thus we
obtain a proof of $\ol\fA(\ft_1) \to \frak{B}$ with only the critical
formulas for $\frak{t}_2$, \dots, $\frak{t}_n$.

Next, replace every formula in the original proof by the conditional
$\fA(\ft_1) \to \frak{F}$, and also replace $\epsilon_a \fA(a)$
everywhere by~$\ft_1$. The initial formulas of the resulting derivation
(except those resulting from critical formulas) are again derivable as
before.  The formulas corresponding to the critical formulas are all
of the form 
\[
\fA(\ft_1) \to (\fA(\ft_i) \to \fA(\ft_1))
\]
which are propositional axioms.  We therefore now have a proof of
$\fA(\ft_1) \to 0 \neq 0$ without critical formulas.  Putting the two
proofs together and applying the law of excluded middle, we have found
a proof of $0 \neq 0$ using only critical formulas for $\ft_2$, \dots,
$\ft_n$.  By induction on $n$, there is a proof of $0 \neq 0$ using no
critical formulas at all. In the resulting proof, we can replace
$\epsilon_x \fA(x)$ everywhere by $0$.\note{This is essentially the
same proof as the one presented as the `Hilbertsche Ansatz' by
\citeN[21]{HilbertBernays:39}. There, the proof is carried out for end
formulas $\fB$ instead of the specific $0 \neq 0$. The only other
difference is that instead of using induction on~$n$, Bernays
constructs one derivation of $\ol\fA(\ft_1) \land \ldots \land
\ol\fA(\ft_n) \to \frak{F}$ and $n$ derivations of $\fA(\ft_i) \to
\frak{F}$, and then applies $n$-fold case distinction.}

In a letter to Ackermann dated October 16, 1929, Bernays discusses this
proof and suggests ways of extending the result to overcome problems
that apparently had led Hilbert to abandon the idea in favour of
consistency proofs using the $\epsilon$-substitution method. The
letter begins with a review of the problems the original idea suffered
from:
\begin{quotation}
While working on the \emph{Grundlagenbuch}, I found myself motivated
to re-think Hilbert's second consistency proof for the
$\epsilon$-axiom, the so-called `failed' proof, and it now
seems to me that it can be fixed after all.

Since I know that it is very easy to overlook something in the area of
proofs like this, I would like to submit my considerations to you for
verification.

The stumbling blocks for the completion of the proof were threefold:
\begin{enumerate}
\item It could happen that due to the replacements needed for the
treatment of one critical formula, a different critical formula lost
its characteristic form without, however, thus resulting in a
derivable formula.
\item Incorporating the second identity axiom, which can be replaced
by the axiom
\[
(G) \qquad a = b \to (\epsilon_x A(x, a) = \epsilon_x A(x, b))
\]
in its application to the $\epsilon$-function [footnote: except in the
harmless application consisting in the substitution of an
$\epsilon$-functional for an individual variable in the identity
axiom]---only $\epsilon_x \fA(x)$ are involved here, where $x$ is an
\emph{individual} variable---caused problems.
\item Sometimes a new $\epsilon$-functional appeared after successful
elimination of an $\epsilon$-functional, so that overall no reduction
was achieved.\note{`\sg Anl"asslich der Arbeit f"ur das
Grundlagenbuch sah ich mich dazu angetrieben, den zweiten Hilbertschen
Wf.-Beweis f"ur das $\epsilon$-Axiom, den sogenannten
"`verungl"uckten"' Beweis, nochmals zu "uberlegen, und es scheint mir
jetzt, dass dieser sich doch richtig stellen l"asst.

Da ich weiss, dass man sich im Gebiete dieser Beweise "ausserst leicht
versieht, so m"ochte ich Ihnen meine "Uberlegung zur Pr"ufung
vorlegen.

Die bisherigen Hindernisse f"ur die Durchf"uhrung des Beweises
bestanden in dreierlei:
\begin{enumerate}
\item Es konnte vorkommen, dass durch die Ersetzungen, die bei der
Behandlung einer kritischen Formel auszuf"uhren waren, eine andere
kritische Formel ihre characteristische Gestalt verlor, ohne doch in
eine beweisbare Formel "uberzugehen.
\item Die Ber"ucksichtigung des zweiten Gleichheits-Axioms, das ja in
seiner Anwendung auf die $\epsilon$-Funktion [Footnote: abgesehen von
der harmlosen Anwendung, bestehend in d. Einsetzung eines
$\epsilon$-Funktionals f"ur eine Grundvariable im Gleichheits-Axiom.]
--- es handelt sich hier immer nur um $\epsilon_a \fA(x)$, wobei $x$
eine \emph{Grund}variable ist---durch das Axiom
\[
(G) \qquad a = b \to (\epsilon_x A(x, a) = \epsilon_x A(x, b))
\]
vertreten werden kann, machte Schwierigkeiten.
\item Es kam vor, dass nach gelungener Ausschaltung eines
$\epsilon$-Funktionals ein anderes $\epsilon$-Funktional hinzurat,
sodass im ganzen keine Reduktion nachweisbar war.\eg'
\end{enumerate}
Bernays to
Ackermann, October 16, 1929. Manuscript, 13 pages. In the possession of
Hans Richard Ackermann. See also \citeN{Ackermann:83}.}
\end{enumerate}
\end{quotation}
The difficulties listed by Bernays arise already for the
$\epsilon$-theorem in the general case; dealing with number theory,
i.e., the induction axiom, in the way outlined requires even further
extensions of the method.  Bernays acknowledges this in the letter,
writing, `With the addition of complete induction the method is no
longer, i.e., at least not immediately, applicable.  For that, your
[Ackermann's] method of total substitution [i.e., a solving
$\epsilon$-substitution] would be the simplest way'.  However, even
if an extension to arithmetic is not immediately available, it seems
that Bernays considered the `second proof' valuable and interesting
enough to fix.  To summarize, there are two difficulties: The
first is that the possibilities in which $\epsilon$-terms can be
nested in one another and in which cross-binding of variables can
occur give rise to difficulties in their elimination.  On the one
hand, we replace $\epsilon_x \fA(x)$ by $\frak{t}_1$ in the
second step.  If $\epsilon$-terms other than $\epsilon_x \fA(x)$, but
which contain $\epsilon_x \fA(x)$, say, $\epsilon_y \fB(y, \epsilon_x
\fA(x))$ are also present, we would obtain from a critical formula
\[
\fB(\frak{s}, \epsilon_x \fA(x)) \to \fB(\epsilon_y \fB(y, \epsilon_x
\fA(x)), \epsilon_x \fA(x))
\]
a formula
\[
\fB(\frak{s}, \frak{t}_1) \to \fB(\epsilon_y \fB(y, \frak{t}_1),
\frak{t}_1)
\]
which is a critical formula for a new $\epsilon$-term (this is
Bernays's point (3)).  On the other hand, the formula $\fA(x)$ might
contain another $\epsilon$-expression, e.g., $\epsilon_y \fB(x, y)$,
in which case the corresponding $\epsilon$-term would be of the form
$\frak{e} \equiv \epsilon_x \fA(x, \epsilon_y \fB(x, y))$.  A critical
formula corresponding to such a term is:
\begin{eqnarray*}
& & \fA(\frak{s}, \epsilon_y \fB(\frak{s}, y)) \to \fA(\epsilon_x \fA(x,
\epsilon_y \fB(x, y)), \epsilon_y \fB(\epsilon_x \fA(x, \epsilon_y
\fB(x, y)), y)), \textrm{ i.e.,} \\
& & \fA(\frak{s}, \epsilon_y \fB(\frak{s}, y)) \to 
\fA(\frak{e}, \epsilon_y \fB(\frak{e}, y))
\end{eqnarray*}
If in this formula the $\epsilon$-term $\epsilon_y \fB(\frak{s}, y)$
is replaced by some other term $\frak{t}$, we get
\begin{eqnarray*}
& & \fA(\frak{s}, \frak{t}) \to \fA(\epsilon_x \fA(x, \epsilon_y \fB(x,
y)), \epsilon_y \fB(\epsilon_x \fA(x, \epsilon_y \fB(x, y)), y)),
\textrm{ i.e.,}\\ & & \fA(\frak{s}, \frak{t}) \to \fA(\frak{e}, \epsilon_y
\fB(\frak{e}, y))
\end{eqnarray*}
which is no longer an instance of the $\epsilon$-axiom. This is
Bernays's point~(1).

The second main difficulty is dealing with equality axioms, for again,
the replacement of an $\epsilon$-term $\epsilon_x \fA(x, a)$ by
$\frak{t}$ might transform an instance of an quality axiom into
\[
a = b \to \frak{t} = \epsilon_x \fA(x, b)
\]
which no longer is an instance of an axiom. (This is Bernays's point~(2)).

Bernays's proposed solution is rather involved and not carried out in
general, but it seems to have prompted Ackermann to apply some methods
from his own \citeyearpar{Ackermann:24} and von Neumann's
\citeyearpar{Neumann:27} $\epsilon$-substitution proofs. Specifically,
the final version of the first $\epsilon$-theorem presented by
\citeN{HilbertBernays:39}, where the solution of the difficulties is
credited to Ackermann, uses double induction on the rank and degree of
$\epsilon$-expressions to deal with the first difficulty, and von
Neumann's notion of $\epsilon$-types to deal with the equality
axiom.\note{The proof of the first epsilon theorem, the `general
consistency result', and Herbrand's theorem of
\citeN{HilbertBernays:39} are also contained in lectures Bernays gave
at Princeton \citeyearpar{Bernays:35a}.}

\section{Hilbert's `Conservativity Programme' and the practice of 
consistency proofs}

A complete understanding of Hilbert's philosophy of mathematics
requires an analysis of what may be called `the practice of finitism'.  
Hilbert was, unfortunately, not always as clear as
one would like in the exposition of his ideas about the finitary
standpoint and of his project of consistency proofs.  Only by
analyzing the approaches by which he and his students attempted to
carry out the consistency programme can we hope to get a complete
picture of the principles and reasoning patterns he accepted as
finitary, and about his views on the nature of logic and proof theory.
The $\epsilon$-substitution method, of course, was considered the most
promising avenue in the quest for a consistency proof.  The `failed
proof' discussed above shows that an alternative approach was, to a
certain degree, pursued in parallel to the more well-known
substitution method, and adds to the understanding we have of Hilbert's
approach to proof theory and consistency proofs.  The `general
consistency result' provides another example of how a consistency
proof should be carried out according to Hilbert.  Its particular
interest lies in its general nature.  Bernays's schematic formulation
of the result underlines and makes explicit the conditions an
axiomatic system should meet in order to be amenable to a consistency
proof of the required form; it stresses once again the requirement of
verifiability and decidability of atomic formulas.

Although Hilbert (and Bernays) have consistently presented the goal of
the proof theoretic programme as giving a finitary consistency proof
of classical, infinitary mathematics, it has become common among
commentators to present the aim of Hilbert's Programme not as one of
finding a proof of consistency, but of finding a proof of
\emph{conservativity of the ideal over the real}.  Here, an appeal is
made to Hilbert's distinction in \citeyearpar{Hilbert:26,Hilbert:28}
between the ideal and real propositions in formalized mathematics: the
real propositions are those that have finitary meaning, whereas
the ideal propositions are those formulas which are added to the real
part to `round out' the theory, to make the uniform application of
logical inferences possible, and which do not have a direct
interpretation on finitist grounds (in particular, they may contain
unbounded quantifiers).  Hilbert likened the real propositions to
those propositions of physical theories which can be verified by
experiment \cite[475]{Hilbert:28}, and hence it is natural to
interpret Hilbert's Programme as an instrumentalist enterprise where
proof theory was supposed to show that whenever a real proposition can
be proved by ideal methods, it can be proved by real methods alone.
Following \citeN{Detlefsen:86}, I shall call this \emph{real-soundness}
of ideal, formalized mathematics.  The first one to present Hilbert's
Programme as aiming for a finitary proof of real-soundness was
\citeN{Neumann:31}. It was probably Kreisel who most consistently and
influentially emphasized it in, e.g.,
\citeyearpar{Kreisel:51,Kreisel:58a,Kreisel:68}.

It must first of all be said that neither Hilbert nor Bernays
presented the aim of the programme as that of finding a proof of
real-soundness; they almost always talk of consistency, and never 
explicitly of
conservativity or of justifying ideal mathematics by showing that
whenever a real statement is provable by ideal means, it is also
finitarily provable.  In fact, there are only two places I could
find where Hilbert formulates something close to a conservativity
claim.  In \citeyearpar{Hilbert:23}, Hilbert writes: `[A] finite theorem
can presumably [\emph{vermutlich}] always be proved without the
transfinite mode of inference [\dots] but this contention is of the
same sort as the contention that every mathematical proposition can
either be proved or refuted'.  The use of `presumably' and the
qualification at the end suggest that Hilbert was not convinced
that the transfinite is conservative over the finite.  Moreover, this
quote appears in the context of a general discussion of mathematical,
not \emph{meta}mathematical proof, and so should not be taken as a
gloss on what the metamathematical consistency proof is supposed to
establish.  The second quote is from \citeN{Hilbert:26}, where he
writes 
\begin{quote}
For there is a condition, a single but absolutely necessary one, to
which the use of the method of ideal elements is subject, and that is
the \emph{proof of consistency}; for, extension by the addition of
ideals is legitimate only if no contradiction is thereby brought about
in the old, narrower domain, that is, if the relations that result for
the old objects whenever the ideal objects are eliminated are valid in
the old domain.  (p.~383)
\end{quote}
This is likewise not as explicit and clear as one would like. Hilbert
is here talking about the method of ideal elements in general, where
this is readily understandable: when we extend real analysis to
complex analysis by extending the domain to include imaginary numbers,
the new theory should not prove any theorems about real numbers which
aren't already true in the real numbers.  The conservativity
requirement formulated here is an explanation of consistency of the
new ideal elements, theorems about them, and their consequences in the
`old, narrower domain', and it is the only place where it is so
formulated.

The question then is: what is the `old, narrower' domain in the case
of proof theory; what are the real sentences for Hilbert?  Hilbert
reserved the label `real formulas' quite clearly for variable-free
formulas.  In \citeyearpar{Hilbert:26}, Hilbert did not use the term
`real' but he did say that the domain being extended by ideal objects
are `formulas to which contentual communications of finitary
propositions [hence, in the main, numerical equations and
inequalities] correspond' (380).  Although general propositions of the
form `for all numerals $\frak{a}$, $\frak{b}$, $\frak{a} + \frak{b} =
\frak{b} + \frak{a}$' are finitary propositions, the corresponding
free-variable formula $a + b = b + a$ is `no longer an immediate
communication of something contentual at all, but a certain formal
object \dots' which does not mean anything in itself (380).  Hilbert's
\citeyearpar{Hilbert:26} was based on a lecture course given in the
Winter semester 1924--25; and in the lecture notes to that course,
Hilbert elaborates on the discussion from which the preceding quote is
taken:
\begin{quotation}
The resulting formulas such as $a + b = b + a$ do not mean anything in
themselves, any more than the numerals do, they are only images of our
thoughts; but from these we can derive propositions, such as $2 + 3 =
3 + 2$ and we are thus led to consider these elementary propositions
also as formulas, and to signify them as such; they then are formulas
which mean those elementary unproblematic propositions.  These
formulas, which mean something, are the old objects, only in a new
conception: all the added formulas, which do not mean anything in
themselves, are the ideal objects of our theory.\note{`\sg Diese so
entstandenen Formeln wie $a + b = b + a$ bedeuten an sich nichts, so
wenig wie die Zahlzeichen, sie sind nur Abbilder unserer Gedanken;
wohl aber k"onnen aus ihnen Aussagen abgeleitet werden, wie $2 + 3 = 3
+ 2$ und wir kommen so dazu, diese elementaren Aussagen auch als
Formeln aufzufassen und zu bezeichnen; es sind das dann Formeln, die
jene elementaren unproblematischen Aussagen bedeuten.  Diese Formeln,
die etwas bedeuten, sind die alten Gebilde, nur in neuer Auffassung;
alle die hinzugef"ugten Formeln hingegen, die an sich nichts bedeuten,
sind die idealen Gebilde unserer Theorie.' \eg\cite[126--127]{Hilbert:24}.}
\end{quotation}
When Hilbert introduces the term `real proposition', he likewise
characterizes them as the `formulas to which correspond contentual
communication of finitary propositions (mainly numerical equations or
inequalities, or more complex communications composed of these)'
\cite[470]{Hilbert:28}.\note{\citeN{Smorynski:89} also holds that the
real propositions are variable-free formulas, and introduces a
three-fold distinction between real propositions (variable-free),
finitary general propositions (quantifier free, but containing free
variables) and ideal propositions (containing quantifiers).  I agree
with \citeN{Detlefsen:90} that Smorynski's reading of Hilbert is
flawed; but do not agree with Detlefsen's assessment of the real/ideal
distinction. This distinction is orthogonal to the finitary/infinitary
distinction.  The real/ideal distinction is a purely syntactic
distinction between those formulas which do not contain variables and
those that do. Every real formula can be immediately interpreted as a
particular finitary propositions, and its truth is decidable. The
ideal formulas cannot immediately be interpreted finitarily, and
when they do \emph{admit} of a finitary interpretation (e.g.,
free-variable formulas as statements about `any given numeral'), these
interpretations are problematic because they do not obey the usual
logical laws (e.g., they are `incapable of being negated').}  Indeed,
if the requirement of real-soundness is understood, as Hilbert and von
Neumann do, by analogy with physical theories and observation
statements, then real propositions must be \emph{decidable}.  So
\citeN[475]{Hilbert:28} says that `only the real propositions are
directly capable of verification', and this only makes sense if `real
proposition' is understood as a variable-free decidable proposition
about numerals.\note{There is no reason to think that `verification'
as used here by Hilbert is the same notion as `\emph{verifizierbar}'
in \cite[238]{HilbertBernays:34}, which does apply to free-variable
formulas and indeed even to formulas with bound variables. See below.}
Since, as Hilbert points out in a slightly different context
\citeyearpar[470]{Hilbert:28}, `one cannot, after all, try out all
numbers' a verification of a general (free-variable) proposition must
consist of a general proof, and as such can hardly be called
`direct'.\note{This interpretation is not undermined by the fact that
people outside the Hilbert school did not always understand `real
proposition' as variable-free; e.g., \citeN{Weyl:28}. In his remarks
at the K\"onigsberg conference \citeyearpar[200--202]{Godel:86},
G\"odel suggests that he understands the real (in his terminology:
`meaningful') propositions as including not only the variable-free
formulas, but also formulas of the form $(Ex)F(x)$.  In his review of
\citeN{Neumann:31}, however, he specifically reports the aim of
Hilbert's Programme as `showing that every numerical formula verifiable
(calculable) in finitely many steps that can be derived according to
the rules of the game by which classical mathematics is played must
turn out to be correct when actually calculated'
\cite[249]{Godel:86}.}

Nevertheless, it has become common among commentators on Hilbert to
take as real all quantifier-free formulas. Since in the formalisms
considered, free-variable formulas are interderivable with their
universal closures, one often takes the real formulas to be all
$\Pi_1$ formulas.  \citeN{Smorynski:77}, for instance, claims that
Hilbert's primary aim was that of establishing conservativity for
$\Pi_1$-formulas, and that he pursued a consistency proof (only?)
because it seemed more tractable and because consistency is equivalent
to conservativity for $\Pi_1$-sentences.  \citeN{Prawitz:81} similarly
states the aim of the programme 
as aiming for `a demonstration in the
real part of mathematics of the fact (if a fact) that every provable
real sentence is true, i.e., that every sentence belonging to the real
part which is proved by possible use of the ideal part is nevertheless
true (according to the standards of the real part)' (254) where `the
real sentences comprise the decidable ones and the ones of the form
$\forall x\, A(x)$ where each instance $A(t)$ is decidable' (256).  He
stated furthermore that this was to be done by proving the claim that
`[f]or each proof $p$ in [ideal mathematics] $T$ and for each real
sentence $A$ in $T$: if $p$ is a proof of $A$ in $T$, then $A$ is
true' (257).  \citeN{Kitcher:76} holds a similar view, and
\citeN[124]{Giaquinto:83} adopts Smory\'nski's formulation of
Hilbert's Programme.  The question of whether $\Pi_1$-conservativity was
a requirement on or an aim of Hilbert's Programme is also of interest
because it underlies an argument against Hilbert's Programme based on
the first incompleteness theorem \cite{Detlefsen:90}. 

As pointed out above, such a formulation of the \emph{aim} of Hilbert's
Programme is not to be found in Hilbert.  Hilbert's sparse remarks, as
well as those contemporary formulations which do focus on
conservativity such as von Neumann's \citeyearpar{Neumann:31}, call
for a proof of conservativity for variable-free, decidable sentences
only. Moreover, and here is where one can bring the historical
discussion in the first part of this paper to bear, the
\emph{practice} of consistency proof only directly establishes
conservativity for quantifier-free sentences.  Consistency proofs
based on the epsilon theorem (e.g., the general consistency result) as
well as those based on $\epsilon$-substitution showed that one can
transform a proof of a closed, variable-free formula (perhaps using
ideal methods, e.g., the transfinite axioms) into a purely real,
variable-free proof which essentially amounts to a calculational
verification of the numerical claim expressed by the end-formula.

It was only Bernays in the \emph{Grundlagen der Mathematik} who drew
the conclusion that the consistency proofs themselves actually
established not only the truth of variable-free formulas provable by
ideal methods, but also of free-variable theorems.  In this context,
Bernays used the term `verifiable' (\emph{verifizierbar}): a
free-variable formula $\fA(a)$ is verifiable if each numerical
instance $\fA(\frak{z})$ is true.  He then stated the following
consequence of consistency proofs: every derivable free-variable
formula is verifiable (This is a consequence of a consistency proof
for quantifier-free formulations of systems of arithmetic in
\citeA[248,298]{HilbertBernays:34}; Bernays also pointed it out as a
consequence of the `general consistency result' in
\citeA[36]{HilbertBernays:39}).  The idea is simple: If $\fA(a)$
(equivalently, $\forall x\, \fA(x)$ is derivable, then the following
method constitutes a finitary proof that, for any $\frak{z}$,
$\fA(\frak{z})$ is true.  From the derivation of $\fA(a)$ we obtain a
derivation of $\fA(\frak{z})$ by substitution. The procedure given in
the consistency proof transforms this derivation into a variable-free
derivation of $\fA(\frak{z})$, which codifies a finitary calculation
that $\fA(\frak{z})$ is true.\note{Proof theorists will realize that
Gentzen's \citeyearpar{Gentzen:36} and Ackermann's
\citeyearpar{Ackermann:40} consistency proofs yield conservativity
result even for $\Pi_2$ sentences, in the follwing sense: If \emph{PA}
proves $\forall x\, \exists y\, A(x, y)$ then there is a
$<\epsilon_0$-recursive function $f$ so that a suitable extension of
\emph{PRA} proves $A(x, f(x))$.  \citeN{Gentzen:36} only vaguely
hinted at this consequence in his paper, saying that the reduction
rules in his consistency proof provide a finitist sense to actualist,
i.e., infinitary propositions.  Again it was
\citeN{Kreisel:51} who expanded on this idea.}

So why has Hilbert been held to $\Pi_1$-conservativity of ideal
mathematics?  \citeN{Kreisel:51} cites Bernays's results; but in
\citeN{Kreisel:58} and later, he instead points to an argument in
\cite[474]{Hilbert:28}. This argument, like Bernays's, shows how a
finitary consistency proof for a system $T$ yields a finitary proof
of every free-variable formula provable in $T$.  Unlike Bernays's
remark, it does not rely on a particular form of the consistency
proof, but on the mere assumption that a finitary consistency proof
is available.  Assume there is a derivation of $\fA(a)$ (equivalently,
of $\forall x\, \fA(x)$).  The task is to show that, for any given
$\frak{z}$, $\fA(\frak{z})$ is true.  Suppose it weren't. Then $\neg
\fA(\frak{z})$ would be true, and, because $T$ proves all true
variable-free formulas, there would be a derivation of $\neg
\fA(\frak{z})$.  But from the derivation of $\fA(a)$ we obtain, by
substitution, a derivation of $\fA(\frak{z})$, and hence $T$ is
inconsistent. But we have a finitary consistency proof of $T$, so
this cannot be the case. Hence, $\fA(\frak{z})$ must be true.  (Note
that this proof uses tertium non datur, but this is a finitarily
admissible application to a variable-free numerical statement
$\fA(\frak{z})$.)  This latter result is presented as a surprising
application of proof theory.  It came several years after the two
models of consistency proofs---$\epsilon$-substitution and the `failed
proof'---had already been worked out. Hence Hilbert here
articulates a conservation property which, \emph{as it turned out,
follows} from consistency---and not, as it were, a property which the
consistency proofs were all along supposed to establish.

To the extent Hilbert saw the original aim of his project as one of
proving conservativity, this aim was to prove conservativity for
decidable, variable-free `real' propositions, but not for
free-variable general propositions.  There is only one passage in
which conservativity for general statements is discussed prior to 1934
(viz., in 1928), and there it is presented as an \emph{application}
and not a statement of the project.  Both the `failed proof' and the
$\epsilon$-substitution method in the writings of Hilbert, Ackermann,
and von Neumann at the time were given the formulation `if there were
an ideal proof of $0 = 1$, then there would be a variable-free (real)
proof of $0=1$'.  If conservativity had played a significant role in
the minds of those involved, it would have been obvious to formulate
the proofs so that they established conservativity.  It is interesting
to note that whereas the $\epsilon$-substitution method only directly
establishes conservativity for variable-free statements (i.e. by
applying it to proofs of variable-free formulas instead of the
specific $0 \neq 1$), the `failed proof' could have been formulated
for proofs of any $\epsilon$-free formula (not just $0 \neq 0$), as it
eventually was in the first $\epsilon$-theorem.  By itself, this would
only have yielded a consistency proof for logic; for arithmetic,
something like the strategy in the `general consistency result' is
needed.  It is nevertheless conceivable that this possible
generalization of the `failed proof' suggested a strategy of how to
remove in general ideal elements in the form of $\epsilon$-terms from
proofs of formulas not involving such ideal terms.  However, the
failed proof was never mentioned in print, and likewise it was never
noted that it would be possible to generalize the then-existing
strategies to give conservativity proofs (in particular,
$\epsilon$-substitution) until well into the 1930s.

In addition to the light they shed on the development of the
conceptual framework of Hilbert's Programme, the results about the
$\epsilon$-calculus discussed above are, I think, of independent and
genuine importance. Interest in the historical development of
Hilbert's Programme has seen a marked increase in the last decade or
so, and naturally the $\epsilon$-calculus takes center stage in the
work on logic in Hilbert's school. Independently of Hilbert studies,
renewed interest in the theory and applications of the
$\epsilon$-calculus \cite{Avigad:02} warrant a closer look at the
foundations and origins of the epsilon calculus---the `failed proof'
is a rather important part of that story.

\section*{Acknowledgments}

Many thanks to Jeremy Avigad and an anonymous referee 
for several helpful suggestions, and to Hans Richard Ackermann 
for permission to quote from the Bernays--Ackermann correspondence.

\theendnotes

\end{document}